\documentclass[a4paper]{article}

\usepackage{amsfonts}
\usepackage{amsmath,amsfonts,amsthm,amssymb}
\usepackage{setspace}
\usepackage{fancyhdr}
\usepackage{lastpage}
\usepackage{extramarks}
\usepackage{chngpage}
\usepackage{tikz}
\usepackage{soul,color}
\usepackage{enumerate}
\usepackage{graphicx,float,wrapfig}

\newtheorem{theorem}{Theorem}

\newtheorem{definition}[theorem]{Definition}
\newtheorem{lemma}[theorem]{Lemma}
\newtheorem{proposition}[theorem]{Proposition}
\newtheorem{corollary}[theorem]{Corollary}
\newtheorem{restate}{Theorem}
\newcommand{\card}[1]{\ensuremath{\left|#1\right|}}

\begin{document}

\title{A New Variation of Hat Guessing Games}

\author{
Tengyu Ma\thanks{Email: matengyu1989@gmail.com} \ \ \ \ \ Xiaoming Sun\thanks{Email: xiaomings@tsinghua.edu.cn} \ \ \ \ \ Huacheng Yu\thanks{Email: yuhch123@gmail.com} \\
\\
Institute for Theoretical Computer Science\\
Tsinghua University, Beijing, China\\
}
\date{}

\maketitle

\begin{abstract}
Several variations of hat guessing games have been popularly discussed in recreational mathematics. In a typical hat guessing game, after initially coordinating a strategy, each of $n$ players is assigned a hat from a given color set. Simultaneously, each player tries to guess the color of his/her own hat by looking at colors of hats worn by other players. In this paper, we consider a new variation of this game, in which we require at least $k$ correct guesses and no wrong guess for the players to win the game, but they can choose to ``pass''.

A strategy is called {\em perfect} if it can achieve the simple upper bound $\frac{n}{n+k}$ of the winning probability. We present sufficient and necessary condition on the parameters $n$ and $k$ for the existence of perfect strategy in the hat guessing games. In fact for any fixed parameter $k$, the existence of perfect strategy can be determined for every sufficiently large $n$.

In our construction we introduce a new notion: $(d_1,d_2)$-regular partition of the boolean hypercube, which is worth to study in its own right. For example, it is related to the $k$-dominating set of the hypercube. It also might be interesting in coding theory. The existence of $(d_1,d_2)$-regular partition is explored in the paper and the existence of perfect $k$-dominating set follows as a corollary.

\end{abstract}

\section{Introduction}

Several different hat guessing games have been studied in recent years \cite{Buhler10, Butler09, Ebert98, Feige04, Feige10, Lenstu05, Maura10}. In this paper we investigate a variation where players can either give a guess or pass. It was first proposed by Todd Ebert in~\cite{Ebert98}. In a standard setting there are $n$ players sitting around a table, who are allowed to coordinate a strategy before the game begins. Each player is assigned a hat whose color is chosen randomly and independently with probability $1/2$ from two possible colors, red and blue. Each player is allowed to see all the hats but his own. Simultaneously, each player guesses its own hat color or passes,  according to their pre-coordinated strategy. If at least one player guesses correctly and no player guesses wrong, the players win the game. Their goal is to design a strategy to maximum their winning probability.

By a simple counting argument there is an upper bound of the maximum winning probability, $n/(n+1)$. It is known that this upper bound can be achieved if and only if $n$ has the form $2^t -1$~\cite{Ebert98}. It turns out that the existence of such perfect strategy that achieves the upper bound corresponds to the existence of perfect 1-bit error-correcting code in $\{0,1\}^n$.

In this paper, we present a natural generalization of Ebert's hat guessing problem: The setting is the same as in the original problem, every player can see all other hats except his own, and is allowed to guess or pass. However, the requirement for them to win the game is generalized to be that at least $k$ players from them should guess correctly, and no player guesses wrong ($1\le k\le n$). Note that when $k=1$, it is exactly the original problem.

We denote by $P_{n,k}$ the maximum winning probability of players. Similarly to the $k=1$ case, $P_{n,k}$ has a simple upper bound $P_{n,k} \le \frac{n}{n+k}$. We call a pair $(n,k)$ {\em perfect} if this upper bound can be achieved, i.e. $P_{n,k}=\frac{n}{n+k}$. There is a simple necessary condition for a pair $(n,k)$ to be perfect, and our main result states that this condition is almost sufficient:

\begin{theorem}\label{IIperf}
For any $d, k, s\in \mathbb{N}$	with $s\geq 2\lceil\lg k\rceil$, $(d(2^s-k), dk)$ is perfect, in particular, $(2^s-k,k)$ is perfect.
\end{theorem}
\noindent There exists pair $(n,k)$ with the necessary condition but not perfect, see the remark in Section~\ref{section:main}.

\vskip0.3cm
Here is the outline of the proof: first we give a general characterization of the winner probability $P_{n,k}$
by using the size of the minimum $k$-dominating set of the hypercube. Then we convert the condition of $(n,k)$ perfect to some kind of {\em regular partition} of the hypercube (see the definition in Section~\ref{section:notation}). Our main contribution is that we present a strong sufficient condition for the existence of such partition, which nearly matches the necessary condition. Then we can transform it into a perfect hat guessing strategy.

As a corollary of Theorem~\ref{IIperf}, we also give asymptotic characterization of the value $P_{n,k}$. For example, we show that for any fixed $k$, the maximum winning probability approaches 1 as $n$ tends to the infinity.

\vskip0.3cm

\noindent {\bf Related work}:

Feige~\cite{Feige04} considered some variations including the discarded hat version and the everywhere balanced version.
 Lenstra and Seroussi~\cite{Lenstu05} studied the case that $n$ is not of form $2^m-1$, they also considered the case with multiple colors.
In \cite{Butler09}, Butler, Hajiaghayi, Kleinberg and Leighton considered the worst case of hat placement with sight graph $G$, in which they need to minimize the maximum wrong guesses over all hat placements.
 In~\cite{Feige10} Feige studied the case that each player can see only some of other players' hats with respect to the sight graph $G$. In \cite{Maura10}, Peterson and Stinson investigated the case that each player can see hats in front of him and they guess one by one.
 Very recently, Buhler, Butler, Graham and Tressler~\cite{Buhler10} studied the case that every player needs to guess and the players win the game if either exactly $k_1$ or $k_2$ players guess correctly, they showed that the simple necessary condition is also sufficient in this game.

\vskip0.3cm
The rest of the paper is organized as follows: Section~\ref{section:notation} describes the definitions, notations and models used in the paper. Then, Section~\ref{section:regular} presents the result of the existence of $(d_1,d_2)$-regular partition of hypercube while Section~\ref{section:main} shows the main result of the hat guessing game. Finally, we conclude the paper in Section~\ref{section:conclusion} with some open problems.

\section{Preliminaries}\label{section:notation}


We use $Q_n$ to denote the the $n$ dimension boolean hypercube $\{0,1\}^n$. Two nodes are adjacent on $Q_n$ if they differ by only one bit. We encode the blue and red color by 0 and 1. Thus the placement of hats on the $n$ players' heads can be represented as a node of $Q_n$. For any $x\in Q_n$,
 $x^{(i)}$ is used to indicate the string obtained by flipping the $i^{th}$ bit of $x$. Throughout the paper, all the operations are over $\mathbb{F}_2$. We will clarify explicitly if ambiguity appears.

Here is the model of the hat guessing game we consider in this paper:  The number of players is denoted by $n$ and players are denoted by $p_1,\ldots,p_n$. The colors of players' hats will be denoted to be $h_1,\dots,h_n$, which are randomly assigned from $\{0,1\}$ with equal probability. $h=(h_1,\ldots,h_n)$. Let $h_{-i} \in Q_{n-1}$ denote the tuple of colors $(h_1,\dots, h_{i-1}, h_{i+1}, \dots, h_n)$ that player $p_i$ sees on the others' heads. The strategy of player $p_i$ is a function $s_i: Q_{n-1}\rightarrow \{0,1, \bot \}$, which maps  the tuple of colors $h_{-i}$ to $p_i$'s answer, where $\bot$ represents $p_i$ answers ``pass" (if some player answers pass, his answer is neither correct nor wrong). A strategy $\mathcal{S}$ is a collection of $n$ functions $(s_1,\ldots,s_n)$. The players win the game if at least $k$ of them guess correctly and no one guesses wrong. We use $P_{n,k}$ to denote the maximum winning probability of the players. The following two definitions are very useful in characterization $P_{n,k}$:

\begin{definition}
A subset $D\subseteq V$ is called a {\em $k$-dominating set} of graph $G=(V,E)$
if for every vertex $v\in V\setminus D$,  $v$ has at least $k$ neighbors
in $D$.
\end{definition}

\begin{definition}
A partition $(V_1,V_2)$ of hypercube $Q_n$ is called a {\em $(d_1,d_2)$-regular partition} if each node in $V_1$ has exactly $d_1$ neighbors in $V_2$, and each node in $V_2$ has exactly $d_2$ neighbors in $V_1$.
\end{definition}

For example, consider the following partition $(V_1,V_2)$ of $Q_3$: $V_1=\{000,111\}$, and $V_2 = Q_3\setminus V_1$. For each vertex in $V_1$, there are $3$ neighbors in $V_2$, and for each vertex in $V_2$, there is exactly one neighbor in $V_1$. Thus $(V_1,V_2)$ forms a $(3,1)$-regular partition of $Q_3$.


\section{($d_1,d_2$)-Regular Partition of $Q_n$} \label{section:regular}

In this section we study the existence of $(d_1,d_2)$-regular partition of $Q_n$.

\begin{proposition}\label{neces_condition}
Suppose $d_1,d_2\le n$, if there exists a $(d_1,d_2)$-regular partition of hypercube $Q_n$, then the parameters $d_1,d_2, n$ should satisfy		$d_1+d_2 = \gcd(d_1, d_2)2^s$ for some $s \le n$.
\end{proposition}

\begin{proof}
Suppose the partition is $(V_1,V_2)$, we count the total number of vertices
	\[
	\card{V_1}+\card{V_2}=2^n,
	\]
	and the number of edges between two parts
	\[
	d_1\card{V_1}=d_2\card{V_2}.
	\]
	By solving the equations, we obtain
	$$\card{V_1}=\frac{d_2}{d_1+d_2}2^n, \ \ \card{V_2}=\frac{d_1}{d_1+d_2}2^n.$$
Both $\card{V_1}$ and $\card{V_2}$ should be integers, therefore $d_1+d_2=\gcd(d_1,d_2)2^s$ holds, since $\gcd(d_1,d_1+d_2)=\gcd(d_2,d_1+d_2)=\gcd(d_1,d_2)$.
\end{proof}

\begin{proposition}\label{l1X}

	If there exists a $(d_1,d_2)$-regular partition of hypercube $Q_n$, then there exists a $(d_1,d_2)$-regular partition of $Q_m$ for every $m \ge n$.
	
\end{proposition}

\begin{proof}
		
	It suffices to show that the statement holds when $m = n+1$, since the desired result follows by induction.  $Q_{n+1}$ can be treated as the union of two copies of $Q_n$ (for example partition according to the last bit), i.e. $Q_{n+1}= Q_{n}^{(1)}\cup Q_{n}^{(2)}$. Suppose $(V_1,V_2)$ is a $(d_1,d_2)$-regular partition of $Q_n^{(1)}$.  We can duplicate the partition $(V_1,V_2)$ to get another partition $(V_1',V_2')$ of $Q_{n}^{(2)}$. Then we can see that $(V_1\cup V_1', V_2\cup V_2')$ forms a partition of $Q_{n+1}$, in which each node has an edge to its duplicate through the last dimension. Observe that each node in $V_1$ ($V_1'$) still has $d_1$ neighbors in $V_2$ ($V_2'$) and same for $V_2$ ($V_2'$), and the new edges introduced by the new dimension are among $V_1$ and $V_1'$, or $V_2$ and $V_2'$, which does not contribute to the edges between two parts of the partition. Therefore we constructed a $(d_1,d_2)$-regular partition of $Q_{n+1}$.
\end{proof}

\begin{proposition}\label{l2X}
	
	If there exists a $(d_1,d_2)$-regular partition of $Q_n$, then there exists $(td_1,td_2)$-regular partition of $Q_{tn}$, for any positive integer $t$.
	
\end{proposition}

\begin{proof}	
	Suppose $(V_1,V_2)$ is a $(d_1,d_2)$-regular partition of $Q_n$. Let $x=x_1x_2\cdots x_{nt}$ be a node in $Q_{nt}$. We can divide $x$ into $n$ sections of length $t$, and denote the sum of $i^{th}$ section by $w_i$, i.e.
	\[
		w_i(x)=\sum_{j=ti-t+1}^{ti} x_j, \ \ (1\leq i\leq n).
	\]
Let $R(x)=w_1(x)w_2(x)\ldots w_n(x)\in Q_n$. Define
$$V_i'=\{x\in Q_{nt}|R(x)\in V_i\}, \ (i=1,2).$$
We claim that $(V_1', V_2')$ is a $(td_1,td_2)$-regular partition of $Q_{nt}$. This is because for any vertex $x$ in $V_1'$, $R(x)$ is in $V_1$. So $R(x)$ has $d_1$ neighbors in $V_2$, and each of which corresponds $t$ neighbors of $x$ in $V_2'$, thus in total $td_1$ neighbors in $V_2'$. It is the same for vertices in $V_2'$.
\end{proof}

By Proposition~\ref{neces_condition}-\ref{l2X} we only need to consider the existence of $(d_1,d_2)$-regular partition of $Q_n$ where $\gcd(d_1,d_2)=1$ and $d_1+d_2=2^s$ (where $s\leq n$), or equivalently, the existence of $(d,2^s-d)$-regular partition of $Q_n$, where $s\leq n$ and $d$ is odd. The following Lemma from~\cite{Buhler10} showed that when $n=2^s-1$ such regular partition always exists.
\begin{lemma}\cite{Buhler10}\label{l3X}
	There exists a $(t,2^s-t)$-regular partition of $Q_{2^s-1}$, for any integer $s,t$ with $0 < t < 2^s$.
\end{lemma}

%
%



The following theorem shows how to construct the $(t,2^s-t)$-regular partition for $n=2^s-r$ (where $r\leq t$).
\begin{theorem}\label{l5X}
	Suppose there exists  a $(t,2^s-t)$-regular partition for $Q_{2^s-r}$ and $t>r$, then there exists a $(t,2^s-t)$-regular partition for $Q_{2^{s+1}-\min\{t,2r\}}$.
\end{theorem}

\begin{proof}
	For convenience, let $m = 2^s-r$, and $l = 2r - \min\{t,2r\}(\ge 0)$. Observe that $2^{s+1} - \min\{t,2r\} = 2m+ l$, and if $t \ge 2r$ then $l = 0$.
	Suppose that $(V_1,V_2)$ is a $(t,2^s-t)$-regular partition for $Q_{m}$.   We want to construct a $(t, 2^{s+1}-t)$-regular partition for $Q_{2m+l}$. The basic idea of the construction is as follows:
	
	We start from set $V_2$. We construct a collection of linear equation systems, each of which corresponds to a node in $V_2$. The variables of the linear systems are the $(2m+l)$ bits of node $x \in Q_{2m+l}$. Let $V_2'$ be the union of solutions of these linear equation systems, and $V_1'$ be the complement of $V_{2}'$. Then $(V_1',V_2')$ is the $(t,2^{s+1}-t)$-regular partition as we desired.
	
\vskip8pt
	Here is the construction. Since $(V_1,V_2)$ is a $(t,2^s-t)$ regular partition for $Q_m$, the subgraph induced by $V_2$ of $Q_m$ is a $(t-r)$-regular graph, i.e. for every node $p\in V_2$, there are $(t-r)$ neighbors of $p$ in $V_2$. For each $p\in V_2$, arbitrarily choose a subset $N(p)\subseteq V_2$ of neighbors of node $p$ with size $\card{N(p)} = r- l$. (here $r-l = r - (2r-\min\{t,2r\}) = \min\{t,2r\} -r$, so $r-l\leq t-r$, and $r-l>0$ since $t > r$)

	Now for each node $p = (p_1,\dots, p_m) \in V_2$, we construct a linear equation system as follows:
	\begin{equation}\label{LinearSystem}
		\begin{cases}
		x_1+x_2=p_1,\\
		x_3+x_4=p_2,\\
		\  \ldots \ \ \ldots \ \ \ldots,\\
		x_{2m-1}+x_{2m}=p_m,\\
		\sum_{j = 1}^{m} x_{2j-1} + \sum_{j \in N(p)} x_{2j} + \sum_{ 1 \le j \le l}x_{2m+j} = 0. \\
		\end{cases}
	\end{equation}
Note that in the last equation the last term $\sum_{ 1 \le j \le l}x_{2m+j}$ vanishes if $l = 0$. Denote by $S(p)\subseteq Q_{2m+l}$ the solutions of this linear system. For convenience, let $f: Q_{2m+l} \rightarrow Q_{m} $ be the operator such that
	\[
		f(x_1, \ldots ,x_{2m+l}) = (x_1+x_2,x_3+x_4,\dots, x_{2m-1}+x_{2m}).
	\]
Then in the linear system~(\ref{LinearSystem}) the first $m$ equations is nothing but $f(x) = p$.
	
	Let $V_2'=\cup_{p\in V_2}S(p)$, and $V_1'=Q_{2m+l}\setminus V_2'$ be its complement. We claim that $(V_1',V_2')$ is a $(t,2^{s+1}-t)$-regular partition of $Q_{2m+l}$.
	
\vskip8pt
	To begin with, observe the following two facts.
	
 	{\bf Observation 1}  For every $x \in V_{2}'$, we have $f(x)  \in V_2$. It can be seen clearly from the first $m$ equations in each equation system.

	{\bf Observation 2} If $k \le 2m$, then $f(x^{(2k)})= f(x^{(2k-1)})  = (f(x))^{(k)} $. If $k>2m$,  $f(x^{(k)}) = f(x)$. Recall the $x^{(i)}$ is the node obtained by flipping the $i^{th}$ bit of $x$.
The observation can be seen from the definition of $f(x)$.

\vskip8pt
	For any node $x \in V_1'$, we show that there are $t$ different ways of flipping a bit of $x$ so that we can get a node in $V_2'$. There are two possible cases:
	
	Case 1: $f(x)  \not \in V_{2}$.  In this case if we flip the $i^{th}$ bit of $x$ for some $i >  2m$, then from Observation 2, $f(x^{(i)})=f(x)$, so $f(x^{(i)})$ will remain not in $V_2$, and therefore $x^{(i)}$ will not be in $V_{2}'$, by Observation 1. So we can only flip the bit in $(x_1,\ldots,x_{2m})$.

	Suppose by flipping the $i^{th}$ bit of $x$ we get $x^{(i)}\in V_2'$ ($i\in [2m]$), from the definition of $V_2'$ we have : $f(x^{(i)}) \in V_2$, and $x^{(i)}$ satisfies the last equation in the equation systems corresponding to $f(x^{(i)})$:
	\begin{equation}\label{LastEquation0}
			\sum_{j = 1}^{m} x_{2j-1} + \sum_{j \in N(f(x^{(i)}))} x_{2j} + \sum_{ 1 \le j \le l}x_{2m+j} = 0.
	\end{equation}
	
	Since $f(x)\notin V_2$ and $(V_1,V_2)$ is a $(t,2^s-t)$-regular partition of $Q_m$, so there are exactly $t$ neighbors of $f(x)$ in $V_2$, which implies there are $t$ bits of $f(x)$ by flipping which we can get a neighbor of $f(x)$ in $V_2$.  Let $\{j_1,\ldots,j_t\}\subseteq [m]$ be these bits, i.e. $f(x)^{(j_1)},\ldots,f(x)^{(j_t)}\in V_2$, by Observation 2,
$$f(x^{(2j_k-1)})=f(x^{(2j_k)})=f(x)^{(j_k)}\in V_2, \ \ (k=1,\ldots,t).$$
But exactly one of $\{x^{(2j_k-1)},x^{(2j_k)}\}$ satisfies the equation~(\ref{LastEquation0}) (here we use the fact $f(x)\notin V_2$, note that $j_k \not \in N(f(x)^{(j_k)})$. Thus totally, there are $t$ possible $i$ such that $x^{(i)} \in V_{2}'$.
	
\vskip8pt
	Case 2: $f(x) \in V_2$. Since $x\notin V_2'$, the last linear equation must be violated, i.e.
 \begin{equation}\label{LastEquation1}
		\sum_{j = 1}^{m} x_{2j-1} + \sum_{j \in N(f(x^{(i)}))} x_{2j} + \sum_{ 1 \le j \le l}x_{2m+j} = 1.
\end{equation}
We further consider three cases here: flip a bit in $\{x_1,\ldots,x_{2m}\}\setminus \{x_{2j},x_{2j-1} : j \in N(f(x))\}$; flip a bit in $ \{x_{2j},x_{2j-1} : j \in N(f(x))\}$; flip a bit in $\{x_{2m+1},\ldots,x_{2m+l}\}$:

 a) if $i \in [m]$ , $i \notin N(f(x))$, and $f(x)^{(i)} \in V_2$. Since $f(x)\in V_2$, $(V_1,V_2)$ is a $(t,2^s-t)$-regular partition of $Q_m$, there are $m-(2^{s}-t)-\card{N(f(x))}=(2^{s} - r) - (2^s-t) - (r-l)= t-2r+l$ such index $i$, and $x^{(2i-1)}$ is the exactly the one in $\{x^{(2i-1)}, x^{(2i)}\}$ which is in $V_{2}'$. (determined by Equation~\ref{LastEquation1}). Thus in this case there are $(t-2r+l)$ neighbors of $x$ in $V_2'$.

 b) if $i \in N(f(x))$, then both of $x^{(2i-1)}, x^{(2i)}$ are in $V_{2}'$, there are $2\cdot \card{N(f(x))} = 2(r-l)$ such neighbors.

 c) if $i > 2m$, then every $x^{(i)}$ is in $V_{2}'$, ($i=2m+1,\ldots,2m+l$), there are $l$ such neighbors.

 Hence, totally $x$ has $(t-2r+l) + 2(r-l) + l = t$  neighbors in $V_2'$.

\vskip8pt
The rest thing is to show that every node $x\in V_2'$ has $(2^{s+1}-t)$ neighbors in $V_1'$. The proof is similar to the proof of Case 2 above, we consider three cases:

a) If $i \in [m]$, $i \not\in N(f(x))$, and $f(x)^{(i)} \in V_2$. Then exactly one of $x^{(2k-1)}, x^{(2k)}$ in $V_{2}'$, thus there are $m-(2^s-t)- \card{N(f(x))}=2^{s} - r - (2^{s}-t) - (r-l)= t-2r+l$ such neighbors of $x$ in $V_{2}'$.

b) If $i \in N(f(x))$ both $x^{(2i-1)}, x^{(2i)}$ are not in $V_{2}'$.

c) If $i >2m$, then every $x^{(i)}$ is not in $V_{2}'$.

Hence totally, $x$ has $(t-2r+l)$ neighbors in $V_{2}'$, and therefore $(2m +l) - (t-2r+l) = 2^{s+1} - t$ neighbors in $V_{1}'$.
	
	Hence we prove that $(V_1',V_2')$ is indeed a $(t,2^{s+1}-t)$-regular partition of $Q_{2^{s+1}-\min\{t,2r\}}$.
\end{proof}

\begin{theorem}\label{l6X}
For any odd number $t$ and any $c \le t$, when $s \geq \lceil \lg t \rceil+\lceil \lg c \rceil $, there exists a $(t,2^s-t)$-regular partition of $Q_{2^s-c}$.
\end{theorem}

\begin{proof}
Let $s_0 =  \lceil \lg t \rceil$. By Lemma~\ref{l3X}, there exists a $(t,2^{s_0}-t)$-regular partition of $Q_{2^{s_0}-1}$. By repeatedly using Theorem~\ref{l5X}, we obtain that there exists $(t,2^{s_0+1}-t)$-regular partition of $Q_{2^{s_0+1}-2}$, $(t,2^{s_0+2}-t)$-regular partition of $Q_{2^{s_0+2}-2^2}$, etc., $(t,2^{s_0+\lceil \lg{c}\rceil-1}-t)$-regular partition of $Q_{2^{s_0+\lceil \lg{c}\rceil-1}-2^{\lceil \lg{c}\rceil-1}}$, and 
$(t,2^{s_0+\lceil \lg{c}\rceil}-t)$-regular partition of $Q_{2^{s_0+\lceil \lg{c}\rceil}-c}$. By using Proposition~\ref{l1X}, we get that there exists a $(t,2^{s}-t)$-regular partition of $Q_{2^{s}-c}$, for any $s \geq s_0 + \lceil \lg{c}\rceil = \lceil \lg t \rceil+\lceil \lg{c}\rceil$.
\end{proof}

Combining Proposition~\ref{l2X} and Theorem~\ref{l6X}, we have the following corollary.
\begin{corollary}\label{regularpartition}
	Suppose $d_1 =dt, d_2 = d(2^s-t), n = d(2^s-c)$, where $d,t,s$ are positive integers with $0 < t < 2^s$, $c \leq t$ and $s \geq \lceil \lg c \rceil + \lceil \lg t \rceil$, then there exists a $(d_1,d_2)$-regular partition for $Q_n$.
\end{corollary}

%
%

\section{The Maximum Winning Probability $P_{n,k}$}
\label{section:main}


The following lemma characterizes the relationship between the maximum winner probability $P_{n,k}$ and the minimum $k$-dominating set of $Q_n$. The same result was showed in~\cite{Feige10} for $k=1$.

\begin{lemma}\label{IIsoldom}
Suppose $D$ is a $k$-dominating set of $Q_n$ with minimum number of vertices. Then
\[P_{n,k} = 1 - \frac{|D|}{2^n}.\]
\end{lemma}
\begin{proof}
Given a $k$-dominating set $D$ of $Q_n$, the following strategy will have winning probability at least $1-\frac{\card{D}}{2^n}$: For any certain placement of hats, each player can see all hats but his own, so player $p_i$ knows that current placement $h$ is one of two adjacent nodes $\{x,x^{(i)}\}$ of $Q_n$. If $x\in D$ (or $x^{(i)}\in D$), he guesses that the current placement is $x^{(i)}$ (or $x$), otherwise he passes.
We claim that by using this strategy, players win the game when the placement is a node which is not in $D$.
Observe that since $D$ is a $k$-dominating set, for any node $y\notin D$, $y$ has $l$ neighbors $y^{(i_1)},y^{(i_2)},\ldots,y^{(i_l)}$ that are in $D$, where $l \ge k$. According to the strategy desribed, players $p_{i_1},\ldots,p_{i_l}$ would guess correctly and all other players will pass. This shows the winning probability is at least $1-\frac{\card{D}}{2^n}$.

\vskip8pt
Next we show that $P_{n,k}\leq 1-\frac{\card{D}}{2^n}$. Suppose we have a strategy with winning probability $P_{n,k}$. We prove that there exists a $k$-dominating set $D_0$, such that $|D_0| = 2^n(1-P_{n,k})$. The construction is straightforward: Let $D_0 = \{ h\in Q_n : h \textrm{ is not a winning placement} \}$. Thus $|D_0| = N(1-P_{n,k})$. For every winning placement $h \notin D_0$, suppose players $p_{i_1},\ldots,p_{i_l}$ will guess correctly ($l\geq k$), consider the placement $h^{(i_1)}$, which differs from $h$ only at player $p_{i_1}$'s hat, so player $p_{i_1}$ will guess incorrectly in this case, thus $h^{(i_1)}\in D_0$. Similarly $h^{(i_2)},\ldots,h^{(i_l)}\in D_0$, therefore $D_0$ is a $k$-dominating set. We have
\[|D| \le |D_0| = 2^n(1-P_{n,k}),\]
which implies
\[
	P_{n,k}\leq 1-\frac{\left|D\right|}{2^n}.
\]

Combining these two results, we have $P_{n,k}=1-\frac{\card{D}}{2^n}$ as desired.
\end{proof}

\begin{proposition}\label{props}
The following properties hold:
\begin{enumerate}[(a)]
	\item
		If $n_1 < n_2$ then $P_{n_1,k} \le P_{n_2,k}$.
	\item
		$(n,k)$ is perfect iff there exists a $(k,n)$-regular partition of $Q_n$.
	\item
		For any $t\in \mathbb{N}$, $P_{nt,kt}\geq P_{n,k}$. As a consequence, if $(n,k)$ is perfect, $(nt,kt)$ is perfect.
\end{enumerate}
\end{proposition}
\begin{proof}
	For part (a), suppose that $D$ is a minimum $k$-dominating set of $Q_{n_1}$. We make $2^{n_2-n_1}$ copies of $Q_{n_1}$, and by combining them we get a $Q_{n_2}$, which has dominating set of size $
2^{n_2-n_1}|D|$. By Lemma \ref{IIsoldom}, $P_{n_2,k} \ge  1
-\frac{2^{n_2-n_1}|D|}{2^{n_2}} = P_{n_1,k}$.

\vskip6pt
	For part (b), 	suppose $(U,V)$ is a $(k,n)$-regular partition of $Q_n$, note that $V$ is a $k$-dominating set of $Q_n$ and $|V| = \frac{k}{n+k}\cdot 2^n$, thus $V$ is a minimum $k$-dominating set of $Q_n$. We have that $P_{n,k} = 1 - \frac{|V|}{2^n} = \frac{n}{n+k}$, which implies that $(n,k)$ is perfect.
	
	On the other hand, if $(n,k)$ is perfect, suppose $D$ is the minimum $k$-dominating set, we have $\left|D\right|=\frac{k}{n+k}\cdot 2^n$. It can be observed that $\left(Q_n\setminus D, D\right)$ is a $(k,n)$-regular partition of $Q_n$.

\vskip6pt	

	For part (c), since $\frac{n}{n+k}=\frac{nt}{nt+kt}$, once $P_{nt,kt}\geq P_{n,k}$ holds, it's an immediate consequence that the perfectness of $(n,k)$ implies the perfectness of $(nk,nt)$.
 	
 	Suppose for $n$ players, we have a strategy $\mathcal{S}$ with probability of winning $P_{n,k}$. For $nt$ players, we divide them into $n$ groups, each of which has $t$ players. Each placement $h=(h_1,h_2,\ldots,h_{nt})$ of $nt$ players can be mapped to a placement $P(h)$ of $n$ players in the following way: for Group $i$, suppose the sum of colors in the group is $w_i$, i.e.
 	\[
 		w_i(h)=\sum_{j=ti-t+1}^{ti} h_j, \ \ (1\leq i\leq n).
 	\]

 	Let $P(h)=(w_1(h),w_2(h),\ldots,w_n(h))$ be a placement of $n$ players. Each player in Group $i$ knows the color of all players in $P(h)$ other than Player $i$, thus he uses Player $i$'s strategy $s_i$ in $\mathcal{S}$ to guess the sum of colors in Group $i$ or passes. Moreover once he knows the sum, his color can be uniquely determined.
 	
 	Note that the players in Group $i$ would guess correctly or incorrectly or pass, if and only if Player $i$ in the $n$-player-game would do. Since the hat placement is uniformly at random, the probability of winning using this strategy is at least $P_{n,k}$, thus $P_{nt,kt}\geq P_{n,k}$.
\end{proof}

Now we can prove our main theorem:

\begin{restate}
For any $d, k, s\in \mathbb{N}$	with $s\geq 2\lceil\lg k\rceil$, $(d(2^s-k), dk)$ is perfect, in particular, $(2^s-k,k)$ is perfect.
\end{restate}

\begin{proof}
	It's an immediate corollary of part (b) of Proposition \ref{props} and Theorem \ref{l6X}.
\end{proof}

\noindent {\bf Remark:}	By Proposition~\ref{neces_condition} and Proposition~\ref{props}(b) there is a simple necessary condition for $(n,k)$ to be perfect, $n+k=\gcd(n,k)2^t$. Theorem \ref{IIperf} indicates that when $n+k=\gcd(n,k)2^t$ and $n$ is sufficiently large, $(n,k)$ is perfect. The necessary condition and sufficient condition nearly match in the sense that for each $k$, there's only a few $n$ that we don't know whether $(n, k)$ is perfect. Moreover, the following proposition shows that the simple necessary condition can't be sufficient. The first counterexample is $(5,3)$, it is not perfect while it satisfies the simple necessary condition. But $(13, 3)$ is perfect by Theorem~\ref{IIperf} and more generally for all $s\geq 4$, $(2^s-3,3)$ is perfect. We verified by computer program that $(2^4-5, 5)=(11,5)$ is not perfect, while by our main theorem $(2^6-5,5)=(59, 5)$ is perfect. But we still don't know whether the case between them, $(2^5-5,5)=(27, 5)$, is perfect.

\begin{proposition}\label{IInotperf}
    $(n,k)$ is not perfect unless $2k+1\leq n$ when $n\geq 2$ and $k<n$.
\end{proposition}

\begin{proof}
	Suppose $(n,k)$ is perfect. According to part (b) of Proposition \ref{props}, we can find $(U,V)$, a $(k,n)$-regular partition of $Q_n$. Suppose $x$ is some node in $U$, and $y$ is some neighbor of $x$ which is also in $U$, $y$ has $k$ neighbors in $V$. They all differ from $x$ at exactly $2$ bits and one of them is what $y$ differs from $x$ at, i.e. each of them ``dominates'' $2$ neighbors of $x$, one of them is $y$. So $x$ has totally $k+1$ neighbors ``dominated'' by $k$ of nodes in $V$. Since all nodes in $V$ are pairwise nonadjacent, these $k+1$ nodes must be in $U$. Now we have $k+1$ neighbors of $x$ are in $U$ and $k$ neighbors are in $V$, it has totally $n$ neighbors. We must have $2k+1\leq n$.
\end{proof}

For each odd number $k$, let $s(k)$ be the smallest number such that $(2^{s(k)}-k,k)$ is perfect. We know that $s(k)\in [\lceil \lg{k}\rceil, 2\lceil \lg{k}\rceil]$.
The following proposition indicates that all $s\geq s(k)$, $(2^s-k,k)$ is also perfect.

\begin{proposition}\label{IIperfmon}
	If $(2^s-k, k)$ is perfect, $(2^{s+1}-k,k)$ is perfect.
\end{proposition}

\begin{proof}
	If $(2^s-k,k)$ is perfect, by Proposition~\ref{props}(b) there is a $(k,2^s-k)$-regular partition of $Q_{2^s-k}$. Thus by Proposition \ref{l1X}, we have a $(k,2^s-k)$-regular partition of $Q_{2^s-k+1}$. Combine this partition and  Theorem~\ref{l5X}, we get a $(k,2^{s+1}-k)$-regular partition of $Q_{2^{s+1}-k}$. Therefore $(2^{s+1}-k, k)$ is perfect.
\end{proof}

Using Theorem~\ref{IIperf} we can give a general lower bound for the winning probability $P_{n,k}$. Recall that there's upper bound $P_{n,k}\le 1-\frac{k}{n+k}$.

\begin{lemma}\label{IIapprox}
	$P_{n,k}>1-\frac{2k}{n+k}$, when $n\geq 2^{2\lceil\lg k\rceil}-k$.
\end{lemma}

\begin{proof}
	Let $n'$ be the largest integer of form $2^t-k$ which is no more than $n$. By Theorem \ref{IIperf}, $(n',k)$ is perfect, i.e. $P_{n',k}=1-\frac{k}{n'+k}$. By part (a) of Proposition \ref{props}, $P_{n,k}\geq P_{n',k}$. On the other hand we have $n+k<2^{t+1}$, so we have
	\[
		P_{n,k}\geq 1-\frac{k}{n'+k}=1-\frac{2k}{2^{t+1}}>1-\frac{2k}{n+k}.
	\]

\end{proof}

\begin{corollary}
	For any integer $k>0$, $\lim_{n\rightarrow \infty}P_{n,k}=1$.
\end{corollary}

\section{Conclusion}\label{section:conclusion}

In this paper we investigated the existence of regular partition for boolean hypercube, and its applications in finding perfect strategies of a new hat guessing games. We showed a sufficient condition for $(n,k)$ to be perfect, which nearly matches the necessary condition. Several problems remain open: for example, determine the minimum value of $s(k)$ such that $(2^{s(k)}-k,k)$ is perfect, and determine the exact value of $P_{n,k}$. It is also very interesting to consider the case when there are more than two colors in the game.


%

\bibliographystyle{plain}

\bibliography{HatGuessing}

\end{document}